\documentclass[12pt]{amsart}
\usepackage{amssymb}
\usepackage{amscd}
\usepackage{amsmath}
\usepackage{enumerate}
\usepackage[all]{xy}
\usepackage[pdftex]{graphicx}

\newtheorem{theo}{Theorem}[section]
\newtheorem{lemma}[theo]{Lemma}

\newtheorem{prop}[theo]{Proposition}

\newtheorem{defi}[theo]{Definition}

\newtheorem{remark}[theo]{Remark}

\newcommand{\B}{{\mathbb{B}}}
\newcommand{\C}{{\mathbb{C}}}

\newcommand{\calA}{{\mathcal{A}}}
\newcommand{\calB}{{\mathcal{B}}}

\newcommand{\calE}{{\mathcal{E}}}

\newcommand{\BS}{{\mathbb{S}}}
\newcommand{\D}{{\mathbb{D}}}
\newcommand{\s}{\vspace{0,3cm}}

\def\e{{\epsilon}}

\def\a{{\alpha}}

\begin{document}
\title{L\^E'S POLYHEDRON FOR LINE SINGULARITIES}
\author{Aur\'elio Menegon Neto}

\thanks{Partial support from CNPq (Brazil).}
\date{14-11-2014}
\address{Aur\'elio Menegon Neto: Departamento de Matem\'atica - Universidade Federal da Para\'iba - Brazil.}
\email{aurelio@mat.ufpb.br}

\begin{abstract}
We study the topology of line singularities, which are complex hypersurface germs with non-isolated singularity given by a smooth curve. We describe the degeneration of its Milnor fiber to the singular hypersurface by means of a vanishing polyhedron in the Milnor fiber. As a milestone, we also study the topology of the degeneration of a complex isolated singularity hypersurface under a non-local point of view.
\end{abstract}
\maketitle

%%%%%%%%%%%%%%%%%%%%%%%%%%%%%%%%%
\section*{Introduction}

The idea of studying the critical level of a complex function by looking at the non-critical level is classical, used by many authors like Milnor, Hirzebruch, Brieskorn, Pham and others. This led to the classic Fibration Theorem of Milnor and to the study of the vanishing homology of a singularity.

In the case of an isolated singularity, L\^e Dung Trang refined in \cite{Le3} the idea of va\-nishing homology and proved that there exists a vanishing polyhedron (triangulable topological space) in the Milnor fiber such that the Milnor fiber is a regular neighbourhood of it, and that there is a continuous map from the Milnor fiber to the singular one which restricts to a homeomorphism outside the polyhedron and takes the polyhedron to the singular point. 

It is unlikely that there is a natural extension of that result to holomorphic functions with arbitrary singular locus. In \cite{MeS} J. Seade and the author proved that there is a vanishing polyhedron in the boundary of the Milnor fiber of any complex hypersurface with one-dimensional singular set. This describes how the link of the singularity is obtained from this boundary, whose topology has been studied by many authors (see \cite{Si2}, \cite{MP2}, \cite{NS} and \cite{BM} for instance).

The main goal of this paper is to show that there is a vanishing polyhedron in the sense of \cite{Le3} for an important class of singularities called {\it line singularities}, defined by D. Siersma in \cite{Si}. These are nothing but complex hypersurface singularity germs with singular set $\Sigma$ a smooth complex curve. 

There is a natural correspondence between line singularities and families of isolated singularity hypersurfaces. In fact, any line singularity $f$ can be seen as a family of isolated singularity hypersurfaces depending holomorphically on the space of parameters $\C$, given by the restrictions of $f$ to a generic family of hyperplane sections $H_s$ transversal at each $s \in \Sigma$. On the other hand, any family of isolated singularity hypersurfaces clearly defines a line singularity. 

This relation motivates the study of line singularities, since it can provide a tool for the study of the topology of families of isolated singularity hypersurfaces. For instance, our construction of a vanishing polyhedron for the line singularity $f$ seems to provide an adequate instrument to study the $\mu$-constant problem (see \cite{BG} and \cite{LR} for instance). This is a work in progress.

The interplay between these two objects leads us to the following definition: We say that a positive real number $\e$ is a {\it good Milnor radius} for the line singularity $f$ as above if $\e$ is a uniform Milnor radius for the corresponding family of isolated singularities (see Definition \ref{defi_4.1}). 

Our main theorem is:

\medskip

\noindent
{\bf Theorem 1.} {\it
Let $f: (\C^{n+1},0) \to (\C,0)$ be a line singularity and let the ball $\B_\e$ around $0$ in $\C^{n+1}$ with radius $\e>0$ be a Milnor ball for $f$. Then there exist:
\begin{itemize}
\item[$(i)$] A polyhedron $\tilde{P}_t$ of real dimension $n+1$ in the Milnor fiber $F_t = f^{-1}(t) \cap \B_\e$, for $t \neq 0$ sufficiently small, such that $F_t$ is a regular neighbourhood of $\tilde{P}_t$;
\item[$(ii)$] A polyhedron $\tilde{P}_0$ in the singular fiber $F_0 = f^{-1}(0) \cap \B_\e$ such that $F_0$ is a regular neighbourhood of $\tilde{P}_0$ and such that either $\tilde{P}_0$ has real dimension $n+1$, if $\e$ is not a good Milnor radius for $f$, or $\tilde{P}_0 = \Sigma \cap \B_\e$ otherwise; 
\item[$(iii)$] A continuous map $\Psi_t: F_t \to F_0$ that takes $\tilde{P}_t$ onto $\tilde{P}_0$ and that restricts to a homeomorphism from $F_t \backslash \tilde{P}_t$ to $F_0 \backslash \tilde{P}_0$;
\item[$(iv)$] A small contractible closed neighbourhood $W$ of $0$ in $\Sigma \cap \B_\e$ and a continuous map $\Upsilon_t: F_t \to (F_0 \cap H_0) \times W$ that takes $\tilde{P}_t$ onto $\{0\} \times W$ and that restricts to a homeomorphism from $F_t \backslash \tilde{P}_t$ to $(F_0 \cap H_0 \backslash \{0\}) \times W$, where $H_0$ is a generic hyperplane section at $0 \in \Sigma$.
\end{itemize}
}

\medskip

We say that a polyhedron $\tilde{P}_t$ as in Theorem 1 is a {\it L\^e polyhedron for $f$} and that a pair of polyhedra $(\tilde{P}_t, \tilde{P}_0)$ as above is a {\it L\^e polyhedral pair for $f$}. As in \cite{Le3}, the map $\Upsilon_t$ gives a geometric realization of the vanishing of the homology $H_*(F_t)$ of $F_t$ to the trivial homology of $F_0$, which gives a geometric realization of the vanishing cycles of $f$. 

To prove this theorem, we will need to consider a not so usual point of view in Singularity Theory, which we call a {\it non-local situation}. It consists on looking at the topology of a holomorphic function $f$ inside a ball that is not necessarily a Milnor ball for $f$, that is, the intersection of $f^{-1}(0)$ with such ball is not necessarily homeomorphic to the cone over its link. This will be the subject of sections \ref{section_1}, \ref{section_2} and \ref{section_3}.

It is expected that Theorem 1 above can be generalized to more general classes of non-isolated singularities. This is a work in progress, in collaboration with J. Seade.

The author is grateful to J. Seade, who introduced him to the subjects concerned in this paper and made significant contributions to it. He also thanks D.T. L\^e and M.A.S. Ruas for many helpful discussions.

Finally, the author thanks the referee for the many constructive comments, corrections and recommendations which helped to improve the readability and quality of the paper.

%%%%%%%%%%%%%%%%%%%%%%%%%%%%%%%%%
\vspace{0.5cm}
\section{The non-local situation for isolated singularities}
\label{section_1}

Let $g: \C^n \to \C$ be a holomorphic function and suppose that $g(0)=0$, in order to simplify notation. It is well known (see \cite{Mi} for instance) that there exists a positive real number $\e>0$ sufficiently small such that for any $\e'$ with $0<\e' \leq \e$ one has that $V(g) := g^{-1}(0)$ intersects the sphere $\BS_{\e'}$ around $0 \in \C^n$ with radius $\e'$ transversally, in the stratified sense. This property gives the so-called {\it conical structure} of the hypersurface $V(g)$ at $0$. A real number $\e>0$ as above is called a {\it Milnor radius} for $g$ at $0$ and the closed ball $\B_\e$ around $0$ with radius $\e$ is said to be a {\it Milnor ball} for $g$ at $0$.

Milnor showed in \cite{Mi} that for any Milnor radius $\e>0$ there exists a positive real number $\eta$, with $0< \eta \ll \e$, such that the restriction:
$$g_|: g^{-1}(\D_{\eta}^*) \cap \B_\e \to \D_{\eta}^*$$
is a locally trivial differentiable fibration, where $\D_{\eta}$ denotes the closed disk around $0$ in $\C$ with radius $\eta$ and $\D_{\eta}^* := \D_{\eta} \backslash \{0\}$. This is the so-called {\it Milnor fibration theorem}. 

Now let $X$ be a reduced $n$-equidimensional complex analytic space, with $0 \in X \subset \C^N$, and let ${\mathcal S} = ({\mathcal S}_\a)_{\a \in A}$ be a Whitney stratification of $X$. In \cite{Le1} L\^e D\~ung Tr\'ang extended the Milnor fibration theorem for any complex analytic function $g: X \to \C$. That is, he showed that if $\e$ is a Milnor radius for $g$ at $0$ and if $g(0) = 0$, then there exists $\eta$, with $0 < \eta \ll \e$, such that the restriction $g_|: g^{-1}(\D_{\eta}^*) \cap \B_\e \to \D_{\eta}^*$ is a locally topologically trivial fibration, where $\B_\e$ is the closed ball around $0$ in $\C^N$ with radius $\e$. This is the so-called {\it L\^e-Milnor fibration theorem}. 

If $\e$ is a Milnor radius for $g$ and $\eta$ is sufficiently small as above, then for any $t \in \D_\eta^*$ the set $X_t := g^{-1}(t) \cap \B_\e$ is called the {\it Milnor fiber of $g$ at $0$}, with boundary $\partial X_t:= X_t \cap \BS_\e$, and the set $X_0 := g^{-1}(0) \cap \B_\e$ is the {\it special fiber of $g$ at $0$}, whose boundary $\partial X_0 := X_0 \cap \BS_\e$ is called the {\it link of $g$ at $0$}. 

For any $t \in \D_\eta$, the Whitney stratification ${\mathcal S}$ of $X$ induces a Whitney stratification of $X_t$ such that $\partial X_t$ is a union of strata. 
%Moreover, since each $X_t$ is analytic we can choose a triangulation of $X_t$ such that $\partial X_t$ is a union of simplices and such that it is compatible with $\mathcal{S}$, that is, each open simplex of $X_t$ is contained in a stratum.
The topological type of $X_t$ does not depend on the Milnor radius $\e$, for any $t \in \D_\eta$ (see Theorem 2.3.1 of \cite{LT}).

We say that $g: X \to \C$ has an isolated singularity at $0 \in X$ if the restriction of $g$ to each stratum ${\mathcal S}_\a$ that does not contain $0$ but whose closure contains $0$ is a submersion and the restriction of $g$ to the stratum $X_{\a(0)}$ that contains $0$ has an isolated critical point at $0$.

D.T. L\^e \cite{Le3} proved the following:

\begin{theo} \label{theo_1.1}
If $g: (X,0) \to (\C,0)$ has an isolated singularity at $0 \in X$ and if $\e$ and $\eta$ are sufficiently small as above, then for each $t \in \D_\eta^*$ there exist:
\begin{itemize}
\item[$(i)$] a polyhedron $P_t$ in $X_t$, compatible with the stratification $\mathcal{S}$, and a continuous map $\tilde{\xi}_t: \partial X_t \to P_t$, compatible with $\mathcal{S}$, such that $X_t$ is homeomorphic to the mapping cylinder of $\tilde{\xi}_t$;
\item[$(ii)$] a continuous map $\Psi_t: X_t \to X_0$ that sends $P_t$ to $\{0\}$ and that restricts to a homeomorphism $X_t \backslash P_t \to X_0 \backslash \{0\}$.
\end{itemize}
Moreover, the construction of the polyhedra $P_t$, the maps $\tilde{\xi}_t$ and the maps $\Psi_t$ can be done simultaneously for all $t$ in a simple path $\gamma \subset \D_\eta$ connecting an arbitrary $t_0 \in \D_\eta^*$ to $0 \in \D_\eta$. This gives a polyhedron $P \subset g^{-1}(\gamma) \cap \B_\e$ such that $g^{-1}(\gamma) \cap \B_\e$ deformation retracts to $P$ and such that $P \cap X_t = P_t$, for any $t \in \gamma \backslash \{0\}$.
\end{theo}

Now we are going to consider the topology of $g$ inside a suitable neighbourhood of $0$ in $\C^N$ that is not necessarily a Milnor ball for $g$. We start defining the relative polar curve of $g$:

\medskip 

For any linear form $\ell: \C^N \to \C$, the restriction of $\ell$ to $X$ induces the analytic morphism:
$$\phi_\ell: X \to \C^2$$
defined by $\phi_\ell(z) = \big( \ell(z), g(z) \big)$, for any $z \in X$. We have:

\begin{lemma} \label{lemma_2.2}
For any compact neighbourhood $W$ of $0$ in $\C^N$, there exists a non-empty Zariski open set $\Omega$ in the space of non-zero linear forms of $\C^N$ to $\C$ that take $0 \in \C^N$ to $0 \in \C$ such that, for any $\ell \in \Omega$, the analytic morphism $\phi_\ell: X \to \C^2$ satisfies: 
\begin{itemize}
\item[$(i)$] The part of the critical locus of the restriction of $\phi_\ell$ to each stratum $S_\a$ that lies in $W \backslash g^{-1}(0)$ is either empty or a complex curve, whose closure we denote by $\Gamma_\a$;
\item[$(ii)$] For each point $p_i$ in the intersection of $\Gamma := \cup_{\a \in A} \Gamma_\a$ with $g^{-1}(0) \cap W$, there exists a small neighbourhood $V_i$ of $p_i$ in $W$ such that $\Gamma \cap V_i$ is a smooth reduced complex curve and such that the restriction of $\phi_\ell$ to $\Gamma \cap V_i$ defines a biholomorphism from $\Gamma \cap V_i$ to its image $\Delta_i:= \phi_\ell(\Gamma \cap V_i)$.
\end{itemize}
\end{lemma}

\begin{proof}
By $(i)$ of Theorem-Definition 1.4.1 of \cite{Le3}, we know that for each $x \in X \cap W$ there exists an open neighbourhood $U_x$ of $x$ in $W$ and a non-empty Zariski open set $\Omega_x$ in the space of non-zero linear forms of $\C^N$ to $\C$ that take $0$ to $0$ such that, for any $\ell \in \Omega_x$, the analytic morphism $\phi_\ell = (\ell,g): X \to \C^2$ is such that the set of the critical points of $\phi_\ell$ that are in $(X \cap U_x) \backslash g^{-1}(0)$ is either empty or a complex curve. Since $W$ is compact, we can choose $x_1, \dots, x_m \in W$ such that $W = U_{x_1} \cup \dots \cup U_{x_m}$. Then the intersection $\Omega := \Omega_{x_1} \cap \dots \cap \Omega_{x_m}$ is a non-empty Zariski open set in the space of non-zero linear forms of $\C^N$ to $\C$ that take $0$ to $0$. Clearly, for any $\ell \in \Omega$, the set of the critical points of $\phi_\ell$ that are in $W \backslash g^{-1}(0)$ is either empty or a complex curve. The proof of $(ii)$ of our lemma follows applying $(ii)$ of Theorem-Definition 1.4.1 of \cite{Le3} at each $p_i$.
\end{proof}

We say that $\ell \in \Omega$ as above is a {\it good linear form for $g$ in $W$} and that $\Gamma_\ell := \cup_{\a \in A} \Gamma_\a$ is the {\it polar curve of $g$ relatively to $\ell$ in $W$}. We also say that $\Delta_\ell := \phi_\ell(\Gamma_\ell)$ is the {\it polar image of $g$ relatively to $\ell$ in $W$}. 

Now we fix some $\ell \in \Omega$. After making a change of coordinates, if necessary, we can assume that $\ell$ is the projection on the first coordinate. We set $\phi := \phi_\ell = (\ell,g)$, \ $\Gamma := \Gamma_\ell$ and $\Delta := \Delta_\ell$, to simplify the notation.

\medskip

We will study the topology of $g$ inside some suitable neighbourhoods of the ambient space $\C^N$, as defined below:

\begin{defi} \label{defi_1.2}
Let $g: X \to \C$ and $\ell: \C^N \to \C$ be as above. We say that a neighbourhood $W_1 \times W_2$ in $\C \times \C^{N-1}$ is an admissible box for $\phi=(\ell,g)$ if:
\begin{enumerate}
\item $W_1$ is a closed ball in $\C$ and $W_2$ is homeomorphic to a closed ball in $\C^{N-1}$;
\item The set of the singular points of $g$ that are in $W_1 \times W_2$ is finite, that is, the restriction of $g$ to each stratum of the Whitney stratification of $X \cap (W_1 \times W_2)$ induced by $\mathcal S$ has only isolated critical points;
\item The boundary of $W_1 \times W_2$ intersects $V(g)$ transversally, in the stratified sense;
\item The hyperplane section $V(g) \cap \ell^{-1}(u)$ intersects $\{u\} \times \partial W_2$ transversally, in the stratified sense, for any $u \in W_1$.
\end{enumerate}
\end{defi}

Although the definition above might seem unnatural at a first glance, one should notice that for any open neighbourhood $U$ in $\C^N$ such that the set of singular points of $g$ that are in $X \cap U$ is finite, there exists an admissible box $W_1 \times W_2$ for $\phi$ contained in $U$. Moreover, if $W_1 \times W_2$ is an admissible box for $\phi=(\ell,g)$, it is also an admissible box for $\phi_s :=(\ell,g_s)$, for any small enough perturbation $g_s$ of $g$. 

\medskip 

Fix an admissible box $W_1 \times W_2$ for $\phi$ and let $\tilde{w}_1, \dots, \tilde{w}_m$ be the singular points of $g$ that are in $W_1 \times W_2$ and set $\{w_1, \dots, w_m\} := \{g(\tilde{w}_1), \dots, g(\tilde{w}_m) \}$. 

Let $\eta>0$ be a small real number such that:

\begin{itemize}
\item[$(A)$] The function $g$ induces a locally topologically trivial fibration:
$$g_|: g^{-1} \big( \D_{\eta} \backslash \{ w_1, \dots, w_m \} \big) \cap (W_1 \times W_2) \to \D_{\eta} \backslash \{ w_1, \dots, w_m \} \, ;$$ 
\item[$(B)$] The map $\phi=(\ell,g)$ induces a locally topologically trivial fibration: 
$$\phi_|: \phi^{-1}(W_1 \times \D_{\eta} \setminus \Delta) \cap (W_1 \times W_2) \to W_1 \times \D_{\eta} \setminus \Delta \, .$$
\end{itemize}

Notice that since $g$ has the Thom $a_f$-property and since $W_1 \times W_2$ satisfies the condition $(3)$ of the definition above, it follows from the first isotopy theorem of Thom-Mather that there exists $\eta_1>0$ such that condition $(A)$ holds. 

Moreover, by condition $(4)$ of Definition \ref{defi_1.2}, we can choose $\eta_2>0$ sufficiently small such that $\phi^{-1}(u,t) = \ell^{-1}(u) \cap g^{-1}(t)$ intersects $\partial (W_1 \times W_2)$ transversally, in the stratified sense, for any $(u,t) \in W_1 \times \D_{\eta_2}$. Then the first isotopy theorem of Thom-Mather implies that condition $(B)$ holds for $\eta_2$.

So the conditions $(A)$ and $(B)$ hold for any $\eta$ with $0<\eta \leq \min\{\eta_1,\eta_2\}$.

\medskip

In the local point of view, one choose $\eta$ so that the intersection $\{ w_1, \dots, w_m \} \cap \D_\eta$ is just the origin. Here we allow $\eta$ to be bigger than this, provided that conditions $(A)$ and $(B)$ above hold.

We remark that in our situation, the topology of $X_t:= g^{-1}(t) \cap (W_1 \times W_2)$ depends on the choice of the neighbourhood $W_1 \times W_2$, for any $t \in \D_\eta$. 

We want to obtain a result like Theorem \ref{theo_1.1} in this case, which we call a {\it non-local situation}. We will prove:

\begin{theo} \label{theo_1.3}
Let $g: X \to \C$ be a holomorphic function and let $W_1 \times W_2$ be an admissible box for $\phi=(\ell,g)$. If $\eta>0$ is small enough as above (i.e. such that conditions $(A)$ and $(B)$ hold), then for each $t \in \D_\eta$ there exist:
\begin{itemize}
\item[$(i)$] A polyhedron $P_t$ in $X_t$, compatible with the stratification $\mathcal{S}$, and a continuous map $\tilde{\xi}_t: \partial X_t \to P_t$, compatible with $\mathcal{S}$, such that $X_t$ is homeomorphic to the mapping cylinder of $\tilde{\xi}_t$;
\item[$(ii)$] A continuous map $\Psi_t: X_t \to X_0$ that sends $P_t$ to $P_0$ and that restricts to a homeomorphism $X_t \backslash P_t \to X_0 \backslash P_0$.
\end{itemize}
Moreover:
\begin{itemize}
\item[$(iii)$] The polyhedron $P_t$ has real dimension $n-1$, for any $t \in \D_\eta \backslash \{ w_1, \dots, w_m \}$.
\item[$(iv)$] Either $P_0 = \{0\}$, if the polar curve $\Gamma$ intersects $g^{-1}(0) \cap (W_1 \times W_2)$ only at $0 \in X$, or $P_0$ has real dimension $n-1$ otherwise.
\item[$(v)$] The construction of the polyhedra $P_t$, the maps $\tilde{\xi}_t$ and the maps $\Psi_t$ can be done simultaneously for all $t$ in a simple path $\gamma \subset \D_\eta$ connecting an arbitrary $t_0 \in \D_\eta \backslash \{ w_1, \dots, w_m \}$ to $0 \in \D_\eta$. This gives a polyhedron $P_\gamma \subset g^{-1}(\gamma) \cap (W_1 \times W_2)$ such that $g^{-1}(\gamma) \cap (W_1 \times W_2)$ is a regular neighbourhood of $P_\gamma$ and such that $P_\gamma \cap X_t = P_t$, for any $t \in \gamma$.
\item[$(vi)$] Actually, the polyhedra $P_t$ and the maps $\tilde{\xi}_t$ can be constructed simultaneously for all $t \in \D_\eta$. This gives a polyhedron $P_\eta \subset g^{-1}(\D_\eta) \cap (W_1 \times W_2)$ such that $g^{-1}(\D_\eta) \cap (W_1 \times W_2)$ is a regular neighbourhood of $P_\eta$ and such that $P_\eta \cap X_t = P_t$, for any $t \in \D_\eta$.
\end{itemize}
\end{theo}

%%%%%%%%%%%%%%%%%%%%%%%%%%%%%%%%%
\vspace{0.5cm}
\section{A special case of Theorem \ref{theo_1.3}}
\label{section_2}

In this section, we will prove a special case of Theorem \ref{theo_1.3}, when $0 \in X$ is the only singularity of $g$ inside $W_1 \times W_2$, that is, when the restriction of $g$ to each stratum of the Whitney stratification of $X \cap (W_1 \times W_2)$ induced by $\mathcal S$ that does not contain $0$ is a submersion and the restriction of $g$ to the stratum that contains $0$ is a submersion outside $0$. 

\medskip

Set: 
$$X_{\D_\eta} := g^{-1}(\D_\eta) \cap (W_1 \times W_2) \, ,$$ 
which equals $\phi^{-1}(W_1 \times \D_{\eta}) \cap (W_1 \times W_2)$ since $\ell(W_1 \times W_2) = W_1$. Notice that $X_t = g^{-1}(t) \cap (W_1 \times W_2)$ equals $\phi^{-1}(W_1 \times \{t\}) \cap (W_1 \times W_2)$, for any $t \in \D_\eta$.

The map $\phi$ induces a map: 
$$\phi_|: X_{\D_\eta} \to W_1 \times \D_{\eta} \, .$$

Suppose that $\eta>0$ is sufficiently small such that conditions $(A)$ and $(B)$ above hold and such that:

\begin{itemize}
\item[$(C)$] The map $\phi$ induces a biholomorphism from $\Gamma \cap X_{\D_\eta}$ to $\phi(\Gamma \cap X_{\D_\eta})$.
\end{itemize}

\begin{figure}[!h] 
\centering 
\includegraphics[scale=0.6]{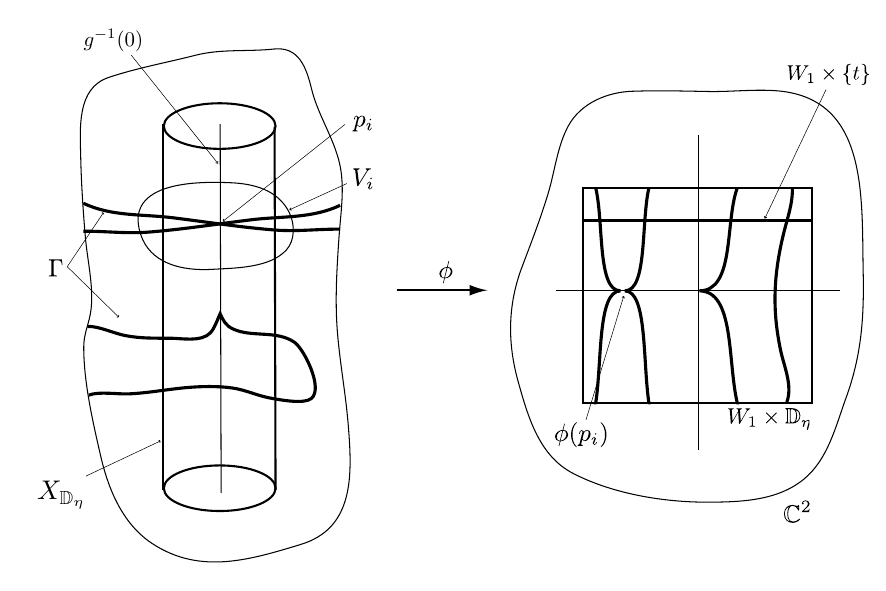}
\caption{}
\label{fig1}
\end{figure}

In fact, to have condition $(C)$ satisfied it is enough to take $\eta>0$ sufficiently small such that $\Gamma \cap X_{\D_\eta}$ is contained in the union $\cup_{i=1}^r V_i$ of the neighbourhoods $V_i$'s given by $(ii)$ of Lemma \ref{lemma_2.2}, putting $W = W_1 \times W_2$. Then it follows from that lemma that $\phi$ induces a biholomorphism from $\Gamma \cap X_{\D_\eta}$ to $\phi(\Gamma \cap X_{\D_\eta})$. See figure \ref{fig1}.

We will prove:

\begin{theo} \label{theo_2.1}
Let $g: X \to \C$ be a holomorphic function with $0 \in X$ and with $g(0)=0$. Let $W_1 \times W_2$ be an admissible box for $\phi=(\ell,g)$ as above. If $\eta>0$ is small enough such that conditions $(A)$, $(B)$ and $(C)$ above hold, then for each $t \in \D_\eta$ there exist:
\begin{itemize}
\item[$(i)$] A polyhedron $P_t$ in $X_t$, compatible with the stratification $\mathcal{S}$, and a continuous map $\tilde{\xi}_t: \partial X_t \to P_t$, compatible with $\mathcal{S}$, such that $X_t$ is homeomorphic to the mapping cylinder of $\tilde{\xi}_t$;
\item[$(ii)$] A continuous map $\Psi_t: X_t \to X_0$ that sends $P_t$ to $P_0$ and that restricts to a homeomorphism $X_t \backslash P_t \to X_0 \backslash P_0$.
\end{itemize}
Moreover:
\begin{itemize}
\item[$(iii)$] The polyhedron $P_t$ has real dimension $n-1$, for any $t \in \D_\eta^*$.
\item[$(iv)$] Either $P_0 = \{0\}$, if the polar curve $\Gamma$ intersects $g^{-1}(0) \cap (W_1 \times W_2)$ in just one point, or $P_0$ has real dimension $n-1$ otherwise.
\item[$(v)$] The construction of the polyhedron $P_t$ and the map $\tilde{\xi}_t$ can be done simultaneously for all $t$ in any closed semi-disk $\D_+ \subset \D_\eta$ containing $0$. This gives a polyhedron $P_+ \subset g^{-1}(\D_+) \cap (W_1 \times W_2)$ such that $g^{-1}(\D_+) \cap (W_1 \times W_2)$ is a regular neighbourhood of $P_+$ and such that $P_+ \cap X_t = P_t$, for any $t \in \D_+$.
\end{itemize}
\end{theo}

The rest of this section is dedicated to the proof of Theorem \ref{theo_2.1}.

\medskip

For any $t \in \D_{\eta}$ set $D_t := W_1 \times \{t\}$ and consider the restriction:
$$\ell_t := \ell|_{X_t}: X_t \to D_t \, .$$
The restriction of $\ell_t$ to each stratum of $X_t$ is a submersion in any point away from $\Gamma$. Therefore it induces a locally trivial fibration over $D_t \backslash (\Delta \cap D_t)$. That is, if we set:
$$\Delta \cap D_t  = \{y_1(t), \dots, y_k(t) \} \, ,$$
then the restriction:
$$\ell_t: X_t \setminus \ell_t^{-1} \big( \{y_1(t), \dots, y_k(t) \} \big) \to D_t \setminus \{y_1(t), \dots, y_k(t) \}$$
is a locally trivial fibration.

\begin{remark}
If $\Gamma$ is empty, one has that $Crit(\phi) \subset g^{-1}(0)$ and then $\phi_|: \phi^{-1} (W_1 \times \D_{\eta}^*) \cap (W_1 \times W_2) \to W_1 \times \D_{\eta}^*$ is a locally topologically trivial fibration, which induces a locally topologically trivial fibration $\ell_t: X_t \to D_t$, for any $t \in \D_{\eta}^*$. Then $X_t$ is homeomorphic to the product of $D_t$ and the general fiber of $\ell_t$. So from now on we shall assume that $\Gamma$ is not empty.
\end{remark}

We will proceed by induction on the dimension of $X$.

%%%%%%%%%%
\s
\subsection{The case when $X$ has dimension $2$.} \label{subsection_dim2}
\ \\

Now we assume that $X$ is a reduced $2$-equidimensional complex analytic space. Then for each $t \in \D_\eta$ fixed, the projection $\ell_t: X_t \to D_t$ induces a finite covering over $D_t \backslash \{ y_1(t), \dots, y_k(t) \}$. 

\begin{figure}[!h] 
\centering 
\includegraphics[scale=0.7]{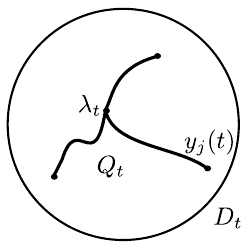}
\caption{}
\label{fig2}
\end{figure}

Let $\lambda_t$ be a point in the interior of $D_t \backslash \{y_1(t), \dots, y_k(t)\}$ and for each $j = 1, \dots, k$, let $\delta(y_j(t))$ be a simple path (differentiable and with no double points) starting at $\lambda_t$ and ending at $y_j(t)$, such that two of them intersect only at $\lambda_t$. See figure \ref{fig2}.

Set $Q_t := \bigcup_{j=1}^k \delta(y_j(t))$ and define:
$$P_t := \ell_t^{-1}(Q_t) \, ,$$
which is a one-dimensional polyhedron in $X_t$, since $\ell_t$ is finite. Clearly, $P_t$ is stratified by the stratification induced by $\mathcal{S}$. 

As in (2.1.1.4) of \cite{Le3}, we can construct a vector field $v_t$ in $D_t$ such that:
\begin{enumerate}
\item It is integrable;
\item It is zero on $Q_t$;
\item It is  transversal to $\partial D_t$ and points inwards;
\item The associated flow $q_t: [0, \infty) \ \times \ (D_t \backslash Q_t) \to D_t$ defines a map
$$
\begin{array}{cccc}
\xi_t \ : & \! \partial D_t & \! \longrightarrow & \! Q_t \\
& \! u & \! \longmapsto & \! \displaystyle \lim_{\tau \to \infty} q_t(\tau,u)
\end{array}
,$$
such that $\xi_t$ is continuous and surjective.
\end{enumerate}

Then we can choose a lifting of  $v_t$ to a vector field $E_t$ in $X_t$ so that: 
\begin{enumerate}
\item It is integrable;
\item It is zero on $P_t$;
\item It is transversal to $\partial X_t$ (in the stratified sense) and points inwards;
\item The associated flow $\tilde{q}_t: [0, \infty) \ \times \ (X_t \backslash P_t) \to X_t$ defines a map
$$
\begin{array}{cccc}
\tilde{\xi}_t \ : & \! \partial X_t & \! \longrightarrow & \! P_t \\
& \! z & \! \longmapsto & \! \displaystyle \lim_{\tau \to \infty} \tilde{q}_t(\tau,z)
\end{array}
,$$
such that $\tilde{\xi}_t$ is continuous, stratified and surjective; 
\item The fiber $X_t$ is homeomorphic to the mapping cylinder of $\tilde{\xi}_t$.
\end{enumerate}

\medskip

In order to construct the collapsing map $\Psi_t$ we do the construction of the vector field $E_t$ simultaneously for all $t$ in a simple path $\gamma$ in $\D_{\eta}$ joining $0$ and some $t_0 \in \partial \D_{\eta}$, such that $\gamma$ is transverse to $\partial \D_{\eta}$. 

The natural projection $\pi: W_1 \times \D_{\eta} \to  \D_{\eta}$ restricted to the polar image $\Delta = \phi(\Gamma)$ induces a ramified covering \index{covering}
$${\pi}_|: \Delta \to \D_{\eta}$$
whose ramification locus is $D_0 \cap \Delta = \{ \phi(p_1), \dots, \phi(p_r) \}$, where $\{p_1, \dots, p_r\}$ is the set of the points in the intersection of $\Gamma$ with $g^{-1}(0) \cap (W_1 \times W_2)$. Notice that in the local situation (when $r=1$) this is just the origin.

\begin{figure}[!h] 
\centering 
\includegraphics[scale=0.6]{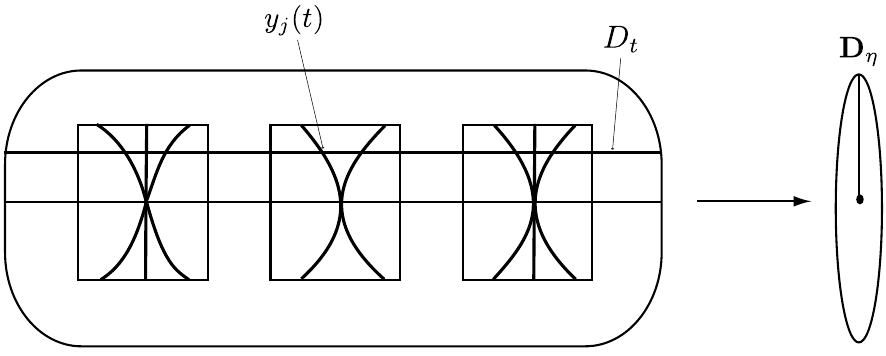}
\caption{}
\label{fig3}
\end{figure}

Hence the inverse image of $\gamma \backslash \{0\}$ by this covering defines $k$ disjoint simple paths in $\Delta$, and each one of them is diffeomorphic to $\gamma \backslash \{0\}$. Each of these paths have $\phi(p_i)$ in its closure, for some $i = 1, \dots, r$, and for any $t \in \gamma \backslash \{0\}$ it contains $y_j(t)$, for some $j = 1, \dots, k$. In particular, we have that $r \leq k$. See figure \ref{fig3}.

We can choose the points $\lambda_t$ above, for each $t \in \gamma$, in such a way that the set $\Lambda := \bigcup_{t \in \gamma} \lambda_t$ defines a simple path in $W_1 \times \D_{\eta}$ such that either $\Lambda \cap \Delta = \{0\}$ if $r=1$, or $\Lambda \cap \Delta = \emptyset$ if $r>1$.

\begin{figure}[!h] 
\centering 
\includegraphics[scale=0.6]{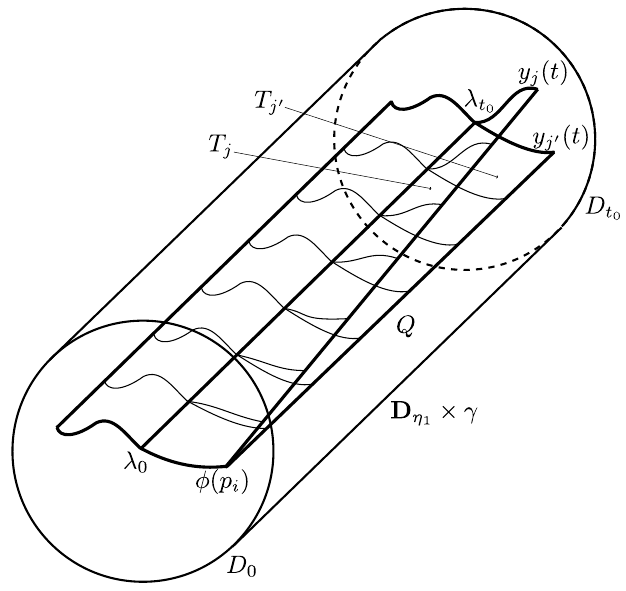}
\caption{}
\label{fig4}
\end{figure}

We can also choose the paths $\delta(y_j(t))$, for each $t \in \gamma$, in such a way that:
$$T_j := \bigcup_{t \in \gamma} \delta(y_j(t))$$
forms either a triangle, if $r=1$, or a square, if $r>1$, immersed in $\bigcup_{t \in \gamma} D_t = W_1 \times \gamma$.  For any $j, j' \in \{1, \dots, k\}$ with $j \neq j'$, note that either $T_j \cap T_{j'} = \Lambda$ or $T_j \cap T_{j'} = \Lambda \cup \gamma(y_j(0)) = \Lambda \cup \gamma(y_{j'}(0))$. See figure \ref{fig4}. 

Set $Q_\gamma := \displaystyle\bigcup_{j=1}^k T_j$ and let $\mathcal{V}$ be a vector field in $W_1 \times \gamma$ such that:
\begin{itemize}
\item[$\bullet$] $\mathcal{V}$ is integrable;
\item[$\bullet$] $\mathcal{V}$ is zero on $Q_\gamma$;
\item[$\bullet$] $\mathcal{V}$ is transversal to $\partial W_1 \times \gamma$; 
\item[$\bullet$] the projection of $\mathcal{V}$ on $\gamma$ is zero.
\end{itemize}

Then the associated flow $w: [0, \infty) \ \times \ (W_1 \times \gamma \backslash Q_\gamma) \to W_1 \times \gamma$ defines a map:
$$ 
\begin{array}{cccc}
\xi \ : & \! \partial W_1 \times \gamma & \! \longrightarrow & \! Q_\gamma \\
& \! z & \! \longmapsto & \! \displaystyle \lim_{\tau \to \infty} w(\tau,z)
\end{array} 
,$$
such that $\xi$ is continuous and surjective.

For any real number $A>0$ set:
$$V_A(Q_\gamma):=  (W_1 \times \gamma) \backslash w \big( [0,A) \times \partial W_1 \times \gamma \big) \, ,$$
which is a closed neighbourhood of $Q_\gamma$ in $W_1 \times \gamma$. Notice that $\partial V_A(Q_\gamma)$ is a differentiable manifold that fibers over $\gamma$ with fiber a circle, and that $W_1 \times \gamma$ is the mapping cylinder of $\xi$. 

Now set: 
$$X_\gamma := \phi^{-1}(W_1 \times \gamma) \cap (W_1 \times W_2) = g^{-1}(\gamma) \cap (W_1 \times W_2) \, .$$
Since the restriction:
$$\phi_|: X_\gamma \backslash \phi^{-1}(Q_\gamma) \to (W_1 \times \gamma) \backslash Q_\gamma$$
is a submersion, it follows that $\phi^{-1}(\partial V_A(Q_\gamma))$ is a submanifold of $X_\gamma$. Notice that the restriction of $\phi$ to $\phi^{-1}(\partial V_A(Q_\gamma))$ is the projection of a locally trivial fibration over $\gamma$. 

Now we set
$$P_\gamma := \phi^{-1}(Q_\gamma) \, .$$
We say that $P_\gamma$ is a {\it collapsing cone for $g$}. Clearly, $P_\gamma$ is a polyhedron of real dimension $2$.

Let $\omega$ be a vector field in $\gamma$ that goes from $t_0$ to $0$ in time $a>0$ and fix $A>0$. We are going to construct an integrable vector field $\calE$ in $X_\gamma \backslash P_\gamma$ such that $\calE$ is tangent to the strata of $\phi^{-1} \big( \partial V_{A'}(Q_\gamma) \big)$, for any $A' \geq A$, in the following way:

\begin{itemize}
\item[$(a)$] If $p \notin \phi^{-1}\big( V_A(Q_\gamma) \big) \cap (W_1 \times W_2)$, there is an open neighbourhood $U_p$ of $p$ in $X_\gamma$ that does not intersect the closed set $\phi^{-1}\big( V_A(Q_\gamma) \big) \cap (W_1 \times W_2)$. We define a differentiable vector field $\calE_p$ on $U_p$ that lifts $\omega$.
\item[$(b)$] If $p \in \left[ \phi^{-1}\big( V_A(Q_\gamma) \big) \cap (W_1 \times W_2) \right] \setminus P_\gamma$, there is an open neighbourhood $U_p$ of $p$ in $X_\gamma$ that does not intersect $P_\gamma$. We define a differentiable vector field $\calE_p$ on $U_p$ that lifts $\omega$ and that is tangent to the strata of $\phi^{-1}\big( \partial V_{A'}(Q_\gamma) \big)$, for any $A' \geq A$.
\end{itemize}
Then the vector field $\calE$ is obtained by gluing the vector fields $\calE_p$, using a partition of unity.

So the flow $h: [0,a] \times X_\gamma \backslash P_\gamma \to X_\gamma \backslash P_\gamma$ associated to $\calE$ defines a homeomorphism $\Psi$ from $X_{t_0} \backslash P_{t_0}$ to $X_0 \backslash P_0$ that extends to a continuous map from $X_{t_0}$ to $X_0$ and that sends $P_{t_0}$ to $P_0$. 

\medskip

Now we can go further and describe the collapsing of $g$ simultaneously, for any $t$ in a closed semi-disk $\D_+$, with $0 \in \D_+ \subset \D_\eta$, in the following way:

For each $j \in \{1, \dots, k\}$ fixed, the union of paths $\delta(y_j(t))$ for all $t \in \D_+$ gives a $3$-dimensional polyhedron $T_j^+$ in $W_1 \times \D_+$. Then we set: 
$$Q_+ := \bigcup_{j=1}^k T_j^+$$
and we consider a vector field $V$ in $W_1 \times \D_+$ that retracts $W_1 \times \D_+$ onto $Q_+$. 

Finally, we set:
$$P_+:= \phi^{-1}(Q_+) \, ,$$
which we call the {\it collapsing polyhedron of $g$ along the semi-disk $\D_+$}. It is a polyhedron of real dimension $3$ contained in $X_+ := \phi^{-1}(W_1 \times \D_+) \cap (W_1 \times W_2)$. 

Since $\phi$ is a submersion over $(W_1 \times \D_+) \backslash Q$, we can lift $V$ to an integrable vector field $E$ in $X_+ \backslash P_+$ that gives the retraction of $X_+$ onto $P_+$.

%%%%%%%%%%%%%%%%%
\s
\subsection{The case when $X$ has dimension $n \geq 3$} \label{subsection_n3}
\ \\

When the complex dimension $n$ of $X$ is greater than $2$, the proof of Theorem \ref{theo_2.1} follows from the next two propositions:

\begin{prop} \label{prop_2.4}
For any $t \in \D_{\eta}$ there exist:
\begin{itemize}
\vspace{0.17cm}
\item[$(i)$] A polyhedron $P_t$ in the fiber $X_t$, adapted to the stratification of $X_t$ induced by $\mathcal{S}$, i.e., the interior of each simplex of $P_t$ is contained in a stratum of $\mathcal{S}$; 
\item[$(ii)$] An integrable vector field $E_t$ in $X_t$, tangent to each stratum of $X_t$, with the properties of the vector field $E_t$ defined in the $2$-dimensional case above.
\end{itemize}
Moreover, $P_t$ has real dimension $n -1$ when $t \neq 0$, and either $P_0$ has real dimension $n -1$ or $P_0 = \{0\}$.
\end{prop}

\medskip

\begin{prop} \label{prop_2.5}
The polyhedron $P_t$ and the vector field $E_t$ of the proposition above can be constructed simultaneously for all $t$ in a closed semi-disk $\D_+$ with $0 \in \D_+ \subset \D_{\eta}$. We obtain a polyhedron $P_+$, adapted to the stratification induced by $\mathcal{S}$, and a continuous vector field $E$ in $X_+ := \phi^{-1}(W_1 \times \D_+) \cap (W_1 \times W_2)$, tangent to each stratum of $X_+$, such that:
\begin{itemize}
\vspace{0.17cm}
\item[$(i)$] For any $t \in \D_+$ the intersection $P_+ \cap X_t$ is a polyhedron $P_t$ and the restriction of $E$ to $X_t$ gives a vector field $E_t$ as in the proposition above;
\item[$(ii)$] The vector field $E$ is tangent to each stratum of $X_+ \backslash P_+$ induced by the stratification $\mathcal{S}$, integrable and non-zero outside $P_+$, zero on $P_+$, transversal to $\partial X_+ := X_+ \cap \partial(W_1 \times W_2)$ in the stratified sense, and pointing inwards.
\end{itemize}
\end{prop}

The propositions above are proved by induction on the complex dimension of $X$: we assume that they are true whenever the dimension of $X$ is $n-1$, and then we prove that they are also true when $X$ has complex dimension $n$.

\medskip

The proof of $(i)$ of Proposition \ref{prop_2.4} follows exactly the same steps as section 3 (pages 312 to 316) of \cite{Le3}. Briefly, we first let $\lambda_0 := (u,0) \in W_1 \times \D_{\eta}$ be the barycenter of $\{ \phi(p_1), \dots, \phi(p_r) \}$ in $W_1 \times \{0\}$. Then we set $\Lambda := \{u\} \times \D_\eta$, which we can suppose either intersects $\Delta$ only at $\lambda_0 = (0,0)$ if $r=1$, or does not intersect $\Delta$ otherwise. Now, given $t_0 \in \D_\eta$ fixed, by induction hypothesis we can consider a L\^e polyhedron $P_{t_0}'$ for the restriction $g'$ of $g$ to the section $X \cap \{\ell=u\}$. If $t_0 \neq 0$, then for each $x_j(t_0) \in \Gamma \cap X_{t_0}$, with $j=1, \dots, k$, we attach a Lefschetz thimble $S_j$ (see \cite{S}, page 221, for instance) that is glued to $P_{t_0}'$ along a sub-polyhedron $(P_j)'_{t_0}$ of $P_{t_0}'$. That is:
$$P_{t_0} = P_{t_0}' \bigcup_{j=1}^k S_j \, .$$
If $t_0=0$ and if $r>1$, then for each $p_i$, with $i=1, \dots, r$, we attach a Lefschetz thimble $S_i$ that is glued to $P_0'$ along a sub-polyhedron $(P_i)'_0$ of $P_0'$. That is:
$$P_0 = P_0' \bigcup_{i=1}^r S_i \, .$$
If $r=1$, then $P_0 = \{0\}$.

\medskip

The proof of $(ii)$ of Proposition \ref{prop_2.4} is identical to the construction of the vector field $E_t$ of \cite{Le3} (subsections 5.1 and 5.2, pages 319 to 326). We will not reproduce it here.

The proof of Theorem \ref{theo_2.1} using the two propositions above is identical to the proof for the $2$-dimensional case. 

\medskip

Now we are going to prove Proposition \ref{prop_2.5}, assuming that Proposition \ref{prop_2.4}, Proposition \ref{prop_2.5} and Theorem \ref{theo_2.1} are true whenever $X$ has complex dimension $n-1$. Actually, we will construct the collapsing polyhedron $P_+$ and we will leave the construction of the vector field $E$ to the reader, since it also follows the steps of subsections 5.1 and 5.2 (pages 319 to 326) of \cite{Le3}. 

\medskip

Fix $t_0 \in \D_+ \backslash \{0\}$ and consider the polyhedron $P_{t_0}$ given by Proposition \ref{prop_2.4}. By the induction hypothesis, we also have a collapsing polyhedron $P_+'$ in $X \cap \{\ell=u\}$ and a continuous integrable vector field $G'$ in $X \cap \{\ell=u\}$ that gives the degeneration of $g'$ along $\D_+$. 

Recall that the map $\phi = (\ell,g)$ induces a locally topologically trivial fibration: 
$$\phi_|: \phi^{-1}(W_1 \times \D_{\eta} \setminus \Delta) \cap (W_1 \times W_2) \to W_1 \times \D_{\eta} \setminus \Delta \, ,$$
where $\Delta$ is the polar image of $g$ relatively to $\ell$ in the admissible box $W_1 \times W_2$. 

Define the $3$-dimensional polyhedra $T_j^+$ in $W_1 \times \D_+$, for each $j=1, \dots, k$, as in subsection \ref{subsection_dim2}. That is:
$$T_j^+ := \bigcup_{t \in \D_+} \delta(y_j(t)) \, ,$$
where $\delta(y_j(t))$ is a simple path connecting $y_j(t)$ and $\lambda_t= (u,t)$.

We are going to construct $P_+$ from the initial polyhedron $P_{t_0}$ as follows: 

\begin{itemize}
\item[$\circ$] If $\Gamma := \cup_{\a \in A} \Gamma_\a$ intersects $g^{-1}(0) \cap (W_1 \times W_2)$ in just one point, then $P_+$ is constructed as in subsection 5.3 (pages 326 and 327) of \cite{Le3}. We briefly describe the construction: 
\ \\
Consider a suitable neighbourhood $U$ of $\phi^{-1}(\Lambda)$ in $\C^N$, conic from $\{0\}$, and for each $x_j(t_0)$ consider a suitable neighbourhood $V_j$ of the component of $\Gamma$ that contains $x_j(t_0)$ in $\C^N$, also conic from $\{0\}$, such that it intersects $U$ away from $\phi^{-1}(\Lambda)$ and such that no two of them intersect.
\ \\
Then we extend $G'$ to $U \cup (\cup_{j} V_j)$ and we let $P_+$ be the union of the orbits of $G'$ that intersect $P_{t_0}$. 

\medskip

\item[$\circ$] If $\Gamma := \cup_{\a \in A} \Gamma_\a$ intersects $g^{-1}(0) \cap (W_1 \times W_2)$ in $r$ points, with $r>1$, then $\delta(y_j(0))$ is a path, for any $j=1, \dots, k$. We fix $j$ and notice that $T_j^+$ is a fiber bundle over $\delta(y_j(0))$. Since $\delta(y_j(0))$ is contractible, we can trivialize $G'$ over it. This gives an integrable vector field $G_j$ in $\phi^{-1}(T_j^+)$. Let $P_+^j$ be the polyhedron in $\phi^{-1}(T_j^+)$ given by the orbits of $G_j$ that intersect $P_{t_0}$. Then we set:
$$P_+ := \bigcup_{j=1}^k P_+^j \, .$$
\end{itemize}

%%%%%%%%%%%%%%%%%%%%%%%%%%%%%%%%%%%%%
\vspace{0.5cm}
\section{The proof of Theorem \ref{theo_1.3}}
\label{section_3}

We will construct a collapsing polyhedron for $g$ in the admissible box $W_1 \times W_2$ along an arbitrary path $\gamma$ in $\D_\eta$ that connects a regular value $t_0$ to the special value $0$, as well as the collapsing polyhedron for $g$ in $W_1 \times W_2$ along all the disk $\D_\eta$. Then the proof of Theorem \ref{theo_1.3} follows from Theorem \ref{theo_2.1} and we leave the details to the reader. 

If the singularities of $g: X \to \C$ lie in the same special fiber $V(g)$, the construction of any L\^e polyhedral pair and of any collapsing polyhedron for $g$ is exactly the same as in the previous section, since the singular points are contained in $\Gamma \cap V(g)$. 

Otherwise, we need to consider the situation when $\D_\eta$ has finitely many special values $0, w_1, \dots, w_m$. In this case, the construction of a L\^e polyhedral pair is still the same, but one should observe that the collapsing along a path $\gamma$ that connects a regular value $t_0$ to the special value $0$ and that passes through special values $w_1, \dots, w_m$ actually describes $2m+1$ degenerations: 
\begin{itemize}

\item[$\circ$] the collapsing of $X_{t_0}$ to $X_{w_1}$ along the sub-path of $\gamma$ connecting $t_0$ and $w_1$;

\item[$\circ$] later, the inverse of the collapsing of $X_{t_2'}$ to $X_{w_1}$, for some regular value $t_2' \in \gamma$ between $w_1$ and $w_2$;

\item[$\circ$] inductively, for each $i=2, \dots, m$ and setting $w_{m+1}:=0$, we have the collapsing of $X_{t_i'}$ to $X_{w_i}$, for some regular value $t_i' \in \gamma$ between $w_{i-1}$ and $w_i$; followed by the inverse of the collapsing of $X_{t_{i+1}'}$ to $X_{w_i}$, for some regular value $t_{i+1}' \in \gamma$ between $w_i$ and $w_{i+1}$;

\item[$\circ$] finally, the collapsing of $X_{t_{m+1}'}$ to $X_0$, for some regular value $t_{m+1}' \in \gamma$ between $w_m$ and $0$. 

\end{itemize}

We say that the corresponding polyhedra $P_{w_1}, \dots, P_{w_m}$ are {\it intermediate polyhedra} of the degeneration of $X_{t_0}$ to $X_0$.

\begin{figure}[!h]
\centering
\includegraphics[scale=0.8]{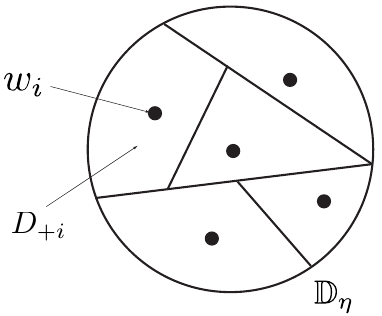}
\caption{}
\label{figure_T5}
\end{figure}

Finally, the construction of the collapsing polyhedron for $g$ along all the disk $\D_\eta$ is done in the following way: Consider a partition of $\D_\eta$ in $m+1$ closed disks $D_{+0}, \dots, D_{+m}$ as in Figure \ref{figure_T5}, in such a way that the intersection of the interior of any two of them is empty and such that $D_{+i}$ contains $w_i$ in its interior, for each $i=0, \dots, m$. 

By the constructions of Theorem \ref{theo_2.1} we get a collapsing polyhedron $P_i$ for $g$ along each $D_{+i}$ such that:

\begin{itemize}
\item The intersection $P_i \cap X_t$ is a L\^e polyhedron $P_t$ for $g$, for any $t \in D_{+i}$ ;
\item If $t \in \partial D_{+i} \cap \partial D_{+j}$ then $P_i \cap X_t = P_j \cap X_t$.
\end{itemize}

So we set $P_\eta:= \bigcup_{i=0}^m P_i$. This completes the proof of Theorem \ref{theo_1.3}.

%%%%%%%%%%%%%%%%%%%%%%%%%%%%%%%%%%%%%
\vspace{0.5cm}
\section{The Proof of Theorem 1}
\label{section_4}

Let $f: (\C^{n+1},0) \to (\C,0)$ be a germ of line singularity. This means that the singular set $\Sigma$ of $V(f)$ is a smooth curve at $0 \in \C^{n+1}$, which we can suppose to be the first coordinate axis $\{z_1 = \dots = z_n=0\}$. 

For each $s \in \C$, let $f_s: \C^n \to \C$ be the restriction of $f$ to the hyperplane section $H_s:= \{z_0=s\}$. It defines a family $(f_s)$ of isolated singularities that depends holomorphically on the complex parameter $s$.

Let $\e$ be a Milnor radius for $f$ and for $f_0$. That is, $\e>0$ is small enough such that $\tilde{\B}_\e \subset \C^{n+1}$ is a Milnor ball for $f$ and such that $\B_\e \subset H_0$ is a Milnor ball for $f_0$. Also let $\ell: \C^n \to \C$ be a generic linear form for $f_0$ as in Section \ref{section_1}, with $\ell(0)=0$. After making another change of coordinates in $\C^{n+1}$, if necessary, we can suppose that $\ell$ is the projection of $H_0$ onto its first coordinate, that is, $\ell(0, z_1, \dots, z_n) = z_1$. 

We can choose a small closed ball with corners $\{0\} \times \D_{\theta_1} \times \D_{\theta_2}$ around $0$ in $H_0$ (which is identified with $\{0\} \times \C \times \C^{n-1}$), contained in the Milnor ball $\B_\e \subset H_0$, such that $\{0\} \times \D_{\theta_1} \times \D_{\theta_2}$ is an admissible box for $\phi_0 := (\ell,f_0)$. Then we can choose a small real number $\omega$, with $0<\omega<\e$, such that $\ell$ is a good linear form for $f_s$ and such that $\{s\} \times \D_{\theta_1} \times \D_{\theta_2}$ is an admissible box for $\phi_s := (\ell,f_s)$, for any $s \in \D_\omega$.

Then we set $T := \D_\omega \times \D_{\theta_1} \times \D_{\theta_2}$, which is a ball with corners around $0$ in $\C^{n+1}$ (identified with $\C \times \C \times \C^{n-1}$). It is well-known that we can suppose $T$ small enough such that the topology of $f$ inside $T$ is equivalent to the topology of $f$ inside the usual Milnor ball $\tilde{\B}_\e \subset \C^{n+1}$ (see Theorem 2.3.1 of \cite{LT} for instance). That is, the fiber of the fibration $f_|: f^{-1}(\D_\eta^*) \cap T \to \D_\eta^*$ is homeomorphic to the fiber of the Milnor fibration $f_|: f^{-1}(\D_\eta^*) \cap \tilde{B}_\e \to \D_\eta^*$, for $\eta$ sufficiently small. 

\begin{defi} \label{defi_4.1}
A positive real number $\e>0$ is a good Milnor radius for $f$ if $\e$ is a Milnor radius for $f$ and if there exists $\omega>0$ sufficiently small such that $\e$ is a Milnor radius for $f_s$, for any $s \in \D_\omega$.
\end{defi}

For any $t_0 \in \D_\eta^*$ fixed, we want to understand the degeneration of $F_{t_0} := f^{-1}(t_0) \cap T$ to the singular fiber $F_0 := f^{-1}(0) \cap T$. 

\begin{figure}[!h]
\centering
\includegraphics[scale=0.8]{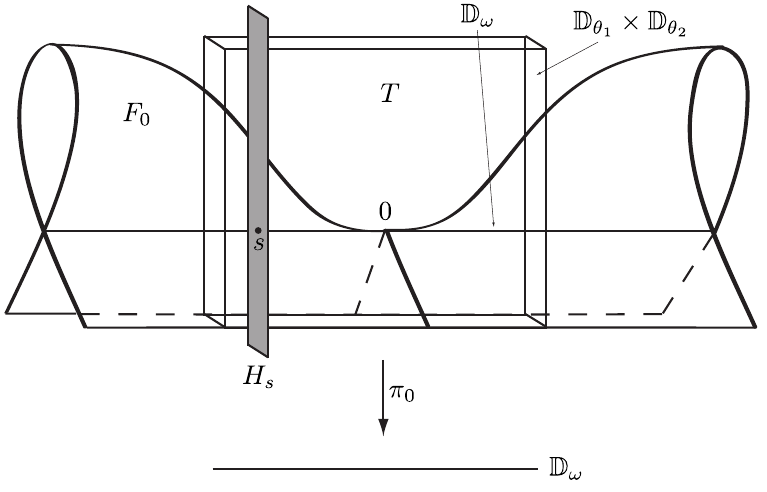}
\caption{}
\label{figure_T3}
\end{figure}

Let: 
$$\pi_{t_0}: F_{t_0}  \to \D_\omega$$
be the restriction of the natural projection $\pi$ of $T$ onto $\D_\omega$. See Figure \ref{figure_T3}.
\medskip

\begin{figure}[!h]
\centering
\includegraphics[scale=0.8]{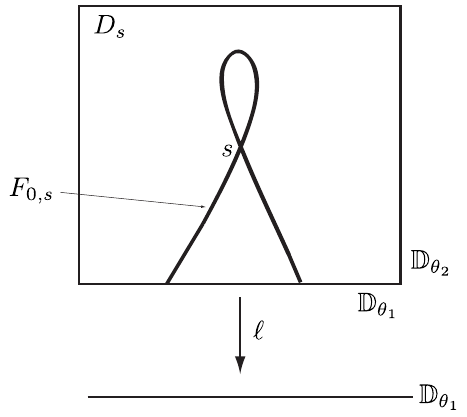}
\caption{}
\label{figure_T4}
\end{figure}

Also, for each $s \in \D_\omega$, set $D_s := T \cap H_s = \{s\} \times \D_{\theta_1} \times \D_{\theta_2}$ and $F_{t_0,s} := \pi_{t_0}^{-1}(s) = f_s^{-1}(t_0) \cap D_s$. See Figure \ref{figure_T4}.

Notice that $F_{t_0,s}$ is not necessarily a Milnor fiber of $f_s$ if $s \neq 0$. Also notice that $\pi_{t_0}$ has isolated critical points in $T$, which are given by the intersection of the polar curve of $f$ relatively to $\pi$ with $F_{t_0}$. Moreover, if we define $\tilde{\ell}: \C^{n+1} \to \C$ by setting $\tilde{\ell} (z_0, z_1, \dots, z_n) = \ell(z_1, \dots, z_n) = z_1$, notice that $\tilde{\ell}(T) = \ell(\D_{\theta_1} \times \D_{\theta_2}) = \D_{\theta_1}$. One can easily check that $T$ is an admissible box for $\phi := (\tilde{\ell}, \pi_{t_0})$. 

Hence we can apply Theorem \ref{theo_1.3} to the projection $\pi_{t_0}$. This gives a L\^e polyhedral pair $(\tilde{P}_{t_0,s}, \tilde{P}_{t_0,0})$ for $\pi_{t_0}$, for any $s \in \D_\omega$, provided we take $\omega>0$ sufficiently small (i.e. $\omega$ is such that conditions $(A)$ and $(B)$ of Section \ref{section_1} hold). Even more, by $(vi)$ of that theorem, we obtain a collapsing polyhedron $\tilde{P}_{t_0} \subset F_{t_0}$ of $\pi_{t_0}$ along $\D_\omega$.  In particular, we have that $F_{t_0}$ is a regular neighbourhood of $\tilde{P}_{t_0}$. We also obtain a vector field $E_{t_0}$ in $F_{t_0}$ that gives the retraction of $F_{t_0}$ onto $\tilde{P}_{t_0}$. Moreover, for any $s \in \D_\omega$ fixed, one has that $\tilde{P}_{t_0,s} = \tilde{P}_{t_0} \cap H_s$ is a L\^e polyhedron for $f_s$, which has real dimension $n-1$. Notice that $\tilde{P}_{t_0}$ has real dimension $n+1$.

The same arguments above work when $t=0$. If $\e$ is not a good Milnor radius for $f$, then $\e$ is not a Milnor radius for $f_s$, for some $s \in \D_\omega^*$ arbitrarily close to $0$, and hence the polar curve of $f_s$ relatively to $\tilde{\ell}$ in $D_s$ intersects $F_0 \cap H_s$ in more than one point. But this intersection coincides with the intersection of the polar curve of $\pi_0$ relatively to $\tilde{\ell}$ in $T$ with $F_0 \cap H_s$. So applying Theorem \ref{theo_1.3} to the projection $\pi_0$ we obtain a collapsing polyhedron $\tilde{P}_0$ in $F_0$ which also has real dimension $n+1$. On the other hand, if $\e$ is a good Milnor radius for $f$, we have that the polar curve of $\pi_0$ relatively to $\tilde{\ell}$ in $T$ is empty, so $\pi_0$ is a locally topologically trivial fibration. Then we set $\tilde{P}_0 := \Sigma \cap T$.

\medskip

We have:

\begin{lemma} \label{lemma_4.3}
Let $\gamma$ be a simple path connecting $t_0 \in \D_\eta^*$ to $0$. There exist:
\begin{itemize}
\item[$(i)$] A polyhedron $P_\gamma$ in $F_\gamma := f^{-1}(\gamma) \cap T$ such that $P_\gamma \cap F_{t_0} = \tilde{P}_{t_0}$ and such that the polyhedron $P_s := P_\gamma \cap H_s$ is a collapsing cone for $f_s$ along $\gamma$, for any $s \in \D_\omega$;
\item[$(ii)$] An integrable vector field $E_\gamma$ in $F_\gamma$ such that for any $s \in \D_\omega$, the restriction of $E_\gamma$ to $H_s$ gives a vector field with the properties of that of Proposition \ref{prop_2.5}.
\end{itemize}
\end{lemma}

\begin{proof}
The idea is to construct the collapsing of each $f_s$ simultaneously, for all $s \in \D_\omega$, starting from the initial polyhedron $\tilde{P}_{t_0}$. 

Set: 
$$\tilde{\phi} := (\pi, \tilde{\ell}, f): \C^{n+1} \to \C^3$$
and:
$$F_\eta:= (\tilde{\phi})^{-1}(\D_\omega \times \D_{\theta_1} \times \D_\eta) \cap T  \, .$$

As before, we have that the Milnor fiber $F_t = f^{-1}(t) \cap \B_\e$ is homeomorphic to $(\tilde{\phi})^{-1}(\D_\omega \times \D_{\theta_1} \times \{t\}) \cap T$. So we reset $F_t := (\tilde{\phi})^{-1}(\D_\omega \times \D_{\theta_1} \times \{t\}) \cap T$ and $F_\gamma := (\tilde{\phi})^{-1}(\D_\omega \times \D_{\theta_1} \times \gamma) \cap T$.

For each $s \in \D_\omega$, let $\Gamma_s$ be the polar curve of $f_s$ relatively to $\ell$ in $D_s = \{s\} \times \D_{\theta_1} \times \D_{\theta_2}$ and let $\Delta_s$ be the corresponding polar image. Set:
$$\Gamma := \bigcup_{s \in \D_\omega} \Gamma_s \subset F_\eta $$
and:
$$\Delta := \bigcup_{s \in \D_\omega} \Delta_s \subset (\D_\omega \times \D_{\theta_1} \times \D_\eta)  \, .$$

We can choose a point $u \in \D_{\theta_1}$ such that the set $\D_\omega \times \{u\} \times \D_\eta$ either intersects $\Delta$ at $\D_\omega \times \{u\} \times \{0\}$, if $\e$ is a good Milnor radius for $f$, or it does not intersect $\Delta$ otherwise.

Now, for each $s \in \D_\omega$ and $j=1, \dots, k_s$ fixed, we can construct the sets $T_{j,s}$ in $\{s\} \times \D_{\theta_1} \times \gamma$ as in section \ref{section_2}, that is, each $T_{j,s} \cap (\{s\} \times \D_{\theta_1} \times \{t\})$, for $t \in \gamma$, is a simple path $\delta(y_{j,s}(t))$ connecting the point $y_{j,s}(t) \in \Delta_s \cap (\{s\} \times \D_{\theta_1} \times \{t\})$ and the point $\lambda_{t,s} := (s, u, t)$. We can do this in such a way that the sets $T_{j,s}$ depend continuously on $s \in \D_\omega$. Set:
$$Q_{\gamma, s} = \bigcup_{j=1}^{k_s} T_{j,s}$$
and:
$$Q_\gamma := \bigcup_{s \in \D_\omega} Q_{\gamma, s}  \, .$$

Now we have all the settings we need to prove Lemma \ref{lemma_4.3} by induction on $n$:

\medskip

Suppose that $n=2$, that is, $f$ is a line singularity defined on $\C^3$. Following the steps of subsection \ref{subsection_dim2}, we get a collapsing cone: 
$$P_\gamma := (\tilde{\phi})^{-1}(Q_\gamma)$$
and a vector field $E_\gamma$ in $F_\gamma$ with the desired properties.

\medskip

Now suppose that the lemma is true whenever $f$ is a line singularity defined on $\C^n$, for some $n \geq 2$ fixed. We will prove it is true for any line singularity $f$ defined on $\C^{n+1}$.

First, consider the polyhedron $\tilde{P}_{t_0} \subset F_{t_0}$ previously constructed. Recall that, for each $s \in \D_\omega$, the polyhedron $P_{t_0,s} = \tilde{P}_{t_0} \cap H_s$ is a L\^e polyhedron for $f_s$.

Let $f'$ be the restriction of $f$ to the hyperplane section $\{ \tilde{\ell} =u \}$. By the induction hypothesis, we have a collapsing polyhedron $P_\gamma'$, with $P_\gamma' \cap F_{t_0} = \tilde{P}_{t_0} \cap \{ \tilde{\ell} = u\}$, and a continuous integrable vector field $G'$ in $\{ \tilde{\ell} =u \}$ that gives the degeneration of $f'$ along $\gamma$. In particular, $P_\gamma' \cap H_s$ and the restriction of $G'$ to $H_s$ give the degeneration of the restriction $f_s'$ of $f_s$ to $\{ \tilde{\ell} =u \} \cap H_s$.

We are going to construct $P_\gamma$ from the initial polyhedron $\tilde{P}_{t_0}$ as follows. For each $s \in \D_\omega$ fixed:

\begin{itemize}
\item[$\circ$] If the polar curve $\Gamma_s$ intersects $f_s^{-1}(0) \cap D_s$ in just one point, we construct $P_{\gamma,s}$ as in \cite{Le3}: 

Set $\Lambda_s := \{s\} \times \{0\} \times \gamma$, which we suppose intersects $\Delta_s$ only at $\tilde{s} := (s,0,0) \in \D_\omega \times \D_{\theta_1} \times \D_\eta$. 

For each $x_{j,s}(t)$ over $y_{j,s}(t)$, with $t \in \gamma$, choose a small radius $r(t)$ such that the set: 
$$\calB_j := \bigcup_{t \in \gamma^*} \B_{r(t)}(x_{j,s}(t))$$
is a neighborhood of $\cup_{t \in \gamma^*}\{ x_{j,s}(t) \}$ in $H_s$, conic from $\tilde{\tilde{s}} := (s,0,0) \in \D_\omega \times \D_{\theta_1} \times \D_{\theta_2}$, for $j=1, \dots, k_s$. To each $\calB_j$ one can associate a neighborhood: 
$$\calA_j := \bigcup_{t \in \gamma} \D_{s(t)}(y_{j,s}(t))$$
in $\{s\} \times \D_{\theta_1} \times \gamma$, conic from $\tilde{s}$. Also let $\mathcal{U}$ be a neighborhood of $\Lambda_s$, conic from $\tilde{\tilde{s}}$, that meets all the $\mathcal{A}_j$'s, but not containing any $y_{j,s}(t)$. See Figure \ref{figure_T2}.

\begin{figure}[!h]
\centering
\includegraphics[scale=0.8]{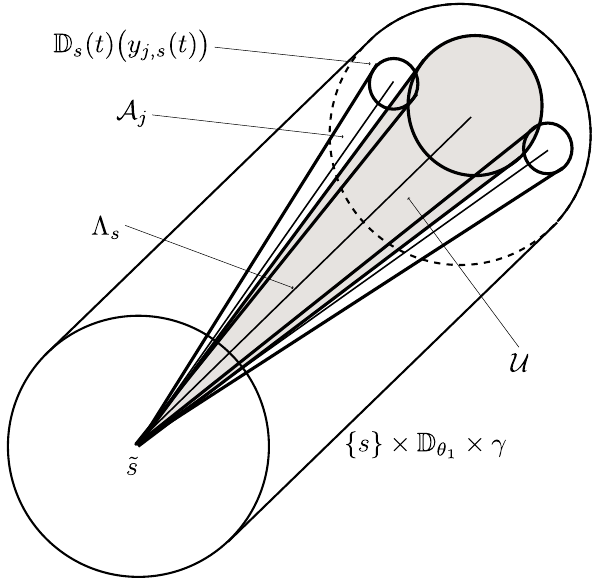}
\caption{}
\label{figure_T2}
\end{figure}

Since $\mathcal{U}$ is contractible, we can extend over $\mathcal{U}$ the vector field $G_s' := G' \cap H_s$. Set:
$$\tilde{\mathcal{U}} := (\tilde{\phi})^{-1}(\mathcal{U}) \cap F_\gamma $$
and let $G_{\tilde{\mathcal{U}}}$ be the vector field in $\tilde{\mathcal{U}}$ given by the extension of $G_s'$. Since $G_s'$ is integrable, the vector field $G_{\tilde{\mathcal{U}}}$ is also integrable.

One can also construct an integrable vector field $G_j$ on each $\calB_j$ that trivializes it over $\gamma$. Then, using a partition of unity, we glue all the vector fields $G_j$'s and $G_{\tilde{\mathcal{U}}}$ together to obtain a trivializing vector field $G_s$ in $\tilde{\mathcal{U}} \cup_{j=1}^k \mathcal{B}_j$ which projects on a radial vector field in $\gamma$ convergent to $0$. This allows us to construct the vanishing cone $P_{\gamma,s}$ from the vanishing polyhedron $P_{t_0,s}$ using the flow of $G_s$.

\medskip

\item[$\circ$] If $\Gamma_s$ intersects $f_s^{-1}(0) \cap D_s$ in more than one point, we construct $P_{\gamma,s}$ as follows: 

Each $y_{j,s}(0)$ is a point of the intersection of $\Gamma_s$ with $f_s^{-1}(0) \cap D_s$, and there are $j_1, j_2 \in \{1, \dots, k_s \}$ with $j_1 \neq j_2$ such that $y_{j_1,s}(0) = y_{j_2,s}(0)$.

Set $\Lambda_s := \{s\} \times \{u\} \times \gamma$, which we suppose intersects $\Delta_s$ only at $\tilde{s} := (s,u,0) \in \D_\omega \times \D_{\theta_1} \times \D_\eta$. Each $\delta(y_{j,s}(0))$ constructed as before is a one-dimensional path, for any $j=1, \dots, k_s$. Fix $j$ and notice that $T_{j,s}$ is a fiber bundle over $\delta(y_{j,s}(0))$. See Figure \ref{fig4}.

Since $\delta(y_{j,s}(0))$ is contractible, we can extend $G'_s$ over it. This gives an integrable vector field $G_{j,s}$ in $(\tilde{\phi})^{-1}(T_{j,s})$. Let $P_{\gamma,s}^j$ be the polyhedron in $(\tilde{\phi})^{-1}(T_{j,s})$ given by the orbits of $G_{j,s}$ that intersect $P_{t_0,s}$. Then we let $P_{\gamma,s}$ be the union of the polyhedra $P_{\gamma,s}^j$, for $j=1, \dots, k_s$.
\end{itemize}

\medskip

Then we set:
$$P_\gamma := \bigcup_{s \in \D_\omega} P_{\gamma,s} \, .$$
Since all the sets considered above depend continuously on $s \in \D_\omega$, we have that $P_\gamma$ is in fact the polyhedron desired.

This allows us to construct the vanishing cone $P_\gamma$ in $F_\gamma$ from a L\^e polyhedron $\tilde{P}_{t_0}$ previously constructed. We leave the construction of the vector field $E_\gamma$ to the reader, since it follows the same arguments as above, together with the steps of the construction of $E_t$ presented in subsections 5.1 and 5.2 of \cite{Le3}. 

\end{proof}

\medskip

Set $\partial_1 F_\gamma := F_\gamma \cap \big( \D_\omega \times \partial (\D_{\theta_1} \times \D_{\theta_2}) \big)$. One can check that the flow $q: [0,\infty) \times (F_\gamma \backslash P_\gamma) \to F_\gamma$ given by the integration of the vector field $E_\gamma$ given by Lemma \ref{lemma_4.3} defines a continuous and surjective map: 
$$ 
\begin{array}{cccc}
\xi \ : & \! \partial_1 F_\gamma & \! \longrightarrow & \! P_\gamma \\
& \! z & \! \longmapsto & \! \displaystyle \lim_{\tau \to \infty} q(\tau,z)
\end{array} 
$$
such that $F_\gamma$ is homeomorphic to the mapping cylinder of $\xi$. This proves $(i)$ and $(ii)$ of Theorem 1.

\medskip

Now we shall construct the collapsing map $\Psi_{t_0}: F_{t_0} \to F_0$. For any positive real number $A>0$ set:
$$V_A(P_\gamma):= F_\gamma \setminus q \big( [0,A) \times \partial_1 F_\gamma \big)$$
which is a closed neighbourhood of $P_\gamma$ in $F_\gamma$. Notice that $\partial V_A(P_\gamma) = q(\{A\} \times \partial_1 F_\gamma)$ is a fiber bundle over $\gamma$, since $\partial_1 F_\gamma$ is a fiber bundle over $\gamma$.

Now let $\sigma$ be a vector field in $\gamma$ that goes from $t_0$ to $0$ in time $a>0$ and fix $A>0$. We are going to construct an integrable vector field $\calE$ on $F_\gamma \backslash P_\gamma$ such that $\calE$ is tangent to $\partial V_{A'}(P_\gamma)$, for any $A' \geq A$, in the following way:

\begin{itemize}
\item[$(a)$] If $p \notin V_A(P_\gamma) \cap T$, there is an open neighbourhood $U_p$ of $p$ in $F_\gamma$ that does not intersect the closed set $V_A(P_\gamma)  \cap T$. So we define a differentiable vector field $\calE_p$ on $U_p$ that lifts $\sigma$.
\item[$(b)$] If $p \in \left[ V_A(P_\gamma) \cap T \right] \setminus P_\gamma$, there is an open neighbourhood $U_p$ of $p$ in $F_\gamma$ that does not intersect $P_\gamma$. Then we can lift $\sigma$ to a differentiable vector field $\calE_p$ on $U_p$ that is tangent to $\partial V_{A'}(P_\gamma)$, for any $A' \geq A$.
\end{itemize}
Then the vector field $\calE$ is obtained by gluing the vector fields $\calE_p$ with a partition of unity.

So the flow $h: [0,a] \times F_\gamma \backslash P_\gamma \to F_\gamma \backslash P_\gamma$ associated to $\calE$ defines a homeomorphism $\Psi_{t_0}$ from $F_{t_0} \backslash \tilde{P}_{t_0}$ to $F_0 \backslash \tilde{P}_0$ that extends to a continuous map from $F_{t_0}$ to $F_0$ and that sends $\tilde{P}_{t_0}$ to $\tilde{P}_0$. This proves $(iii)$ of Theorem 1.

\medskip

Finally, consider the projection $\pi_0: F_0 \to \D_\omega$. Applying Theorem \ref{theo_1.3}, for each $s \in \D_\omega$ we obtain a collapsing map $\Theta_s: F_{0,s} \to F_{0,0}$ that sends $P_{0,s}$ onto $\{0\}$ and that restricts to a homeomorphism $F_{0,s} \backslash P_{0,s} \to F_{0,0} \backslash \{0\}$. Actually, we can do this construction simultaneously for all $s \in \D_\omega$. Hence we have a continuous map $\Theta :  F_0 \to F_{0,0} \times \D_\omega$ defined by:
$$\Theta(z) := (\Theta_{\pi_0(z)} (z) , \pi_0(z))$$
Moreover, $\Theta$ takes $\tilde{P}_0$ onto $\{0\} \times \D_\omega$ and it restricts to a homeomorphism $F_0 \backslash \tilde{P}_0 \to F_{0,0} \backslash \{0\} \times \D_\omega$.

So the map $\Upsilon_t :=  \Theta \circ \Psi_t: F_t \to (F_0 \cap H_0) \times \D_\omega$ is continuous. Moreover, it takes $\tilde{P}_t$ onto $\{0\} \times \D_\omega$ and it restricts to a homeomorphism from $F_t \backslash \tilde{P}_t$ to $F_{0,0} \backslash \{0\} \times \D_\omega$. This proves $(iv)$ of Theorem 1.

%%%%%%%%%%%%%%%%%%%%%%%%%%%%%%%%%

\end{document}